\documentclass[10pt]{amsart}

\usepackage[utf8]{inputenc}

\usepackage{url}

\usepackage{hyperref}

\usepackage{rotating}

\usepackage[small, width=0.9\textwidth]{caption}

\usepackage[all]{xy}
\usepackage{pb-diagram, pb-xy}

\usepackage[dvispnames]{xcolor}

\usepackage{amssymb}
\usepackage{amsthm}

\usepackage{mathrsfs}

\numberwithin{equation}{section}

\theoremstyle{plain}
\newtheorem{thm}[equation]{Theorem}
\newtheorem{lem}[equation]{Lemma}

\newtheorem{prop}[equation]{Proposition}

\newtheorem{rem}[equation]{Remark}

\newtheorem*{thm*}{Theorem}

\theoremstyle{definition}

\newtheorem{ex}[equation]{Example}

\newcommand{\bbC}{{\mathbb C}}

\newcommand{\bbG}{{\mathbb G}}

\newcommand{\bbN}{{\mathbb N}}
\newcommand{\bbP}{{\mathbb P}}

\newcommand{\bbZ}{{\mathbb Z}}

\newcommand{\sfX}{{\mathsf X}}

\newcommand{\rmP}{{\mathrm P}}
\newcommand{\rmS}{{\mathrm S}}

\newcommand{\scrL}{{\mathscr L}}

\newcommand{\frakg}{{\mathfrak g}}

\newcommand{\frakS}{{\mathfrak S}}

\newcommand{\Vect}{\mathrm{Vect}}

\newcommand{\Supp}{\mathrm{Supp}}
\newcommand{\Stab}{\mathrm{Stab}}

\newcommand{\CH}{\mathrm{CH}}

\newcommand{\Alt}{\mathrm{Alt}}

\newcommand{\Lie}{\mathrm{Lie}}

\newcommand{\one}{{\mathbf{1}}}

\newcommand{\Perv}{\mathrm{Perv}}

\newcommand{\Gl}{\mathit{Gl}}
\newcommand{\Sp}{\mathit{Sp}}

\newcommand{\Sl}{\mathit{Sl}}

\newcommand{\Hom}{{\mathrm{Hom}}}

\newcommand{\Aut}{{\mathrm{Aut}}}
\newcommand{\id}{{\mathit{id}}}

\newcommand{\CartierDual}[1]{{{\mathrm{Hom}}(#1, \bbG_m)}}

\newcommand{\scrHom}{{\mathscr{H}\kern -.9pt om}}
\newcommand{\scrExt}{{\mathscr{E}\kern -.9pt xt}}

\renewcommand{\Im}{\mathrm{Im}}

\DeclareMathOperator{\Rep}{{Rep}}

\DeclareMathOperator{\cc}{cc}
\DeclareMathOperator{\cha}{char}
\DeclareMathOperator{\ch}{ch}

\newcommand{\updots}{\text{\reflectbox{$\ddots$}}}

\makeatletter
\def\@seccntformat#1{%
  \protect\textup{\protect\@secnumfont
    \ifnum\pdfstrcmp{subsection}{#1}=0 \bfseries\fi
    \csname the#1\endcsname
    \protect\@secnumpunct
  }%
}  
\makeatother

\begin{document}

\title[]{Summands of theta divisors on Jacobians}
\author{Thomas Kr\"amer}
\address{Institut f\"ur Mathematik, Humboldt-Universit\"at zu Berlin \newline \hspace*{1em} Unter den Linden 6, 10099 Berlin (Germany)}
\email{thomas.kraemer@math.hu-berlin.de}


\begin{abstract}
We show that the only summands of theta divisors on Jacobians of curves and intermediate Jacobians of cubic threefolds are the obvious powers of the curve and the Fano surface of lines on the threefold. The proof uses the decomposition theorem for perverse sheaves, some representation theory and a computation of characteristic cycles for Brill-Noether sheaves.
\end{abstract}

\maketitle
\setcounter{tocdepth}{1}
\tableofcontents

\thispagestyle{empty}

\section{Introduction}

Theta divisors on Jacobian varieties have been studied a lot for their particular properties. If $C$ is a smooth complex projective curve of genus $g\geq 2$ and $A=JC$ denotes its Jacobian variety, then Riemann's theorem says that the theta divisor 
is a sum
\[
 \Theta \;=\; C + \cdots + C \;\subset\; A
\]
of $g-1$ copies of the curve when the latter is embedded into the Jacobian via a suitable translate of the Abel-Jacobi map. More generally, for $d\in\{0,1,\dots, g-1\}$ let $W_d \subset A$ be the image of the sum map $C^d \rightarrow A$ for the above translate of the Abel-Jacobi embedding, then the theta divisor admits the obvious decompositions as a sum
\[
 \Theta 
 \;=\; W_d + W_{g-1-d} 
 \;=\; -W_d - W_{g-1-d}. 
\]
In this note we prove the following converse of Riemann's theorem:

\begin{thm} \label{thm:JC}
If $A$ is the Jacobian of a smooth projective curve of genus $g\geq 2$, then up to translation of the summands the above are the only decompositions of the theta divisor as a sum 
$\Theta = X + Y$ of subvarieties $X,Y\subset A$.
\end{thm}

\pagebreak

This can be seen as a counterpart to Debarre's result on subvarieties of minimal cohomology class in Jacobians~\cite{DebMinimal}. For summands with $\dim X = 1$ it is due to Schreieder who showed that the existence of a curve summand of the theta divisor characterizes Jacobians among all principally polarized abelian varieties~\cite{SchreiederCurveSummands}; this seems to be the only previously known case. There is one other famous example of theta divisors with summands: For any smooth cubic threefold $T\subset \bbP^4$ the intermediate Jacobian~$A=JT$ is a principally polarized abelian fivefold, and it has been shown by Clemens and Griffiths~\cite{CGIntermediateJacobian} that its theta divisor can be decomposed as a sum
\[
 \Theta 
 \;=\; S + (-S) 
 \;=\; (-S) + S
\]
where $S\subset A$ is a copy of the Fano surface of the threefold. Again we show:

\begin{thm} \label{thm:JT}
If $A$ is the intermediate Jacobian of a smooth cubic threefold, then up to translation of the summands the above are the only decompositions of the theta divisor as a sum $\Theta = X+Y$ of positive-dimensional subvarieties $X,Y\subset A$.
\end{thm}

Under the additional assumption that the summands have minimal cohomology class, this is a result of Casalaina-Martin, Popa and Schreieder~\cite{CMPSGenericVanishing}; for a further discussion of summands of theta divisors and generic vanishing subschemes we refer to~\cite{SchreiederDecomposable}. Our approach is very different, it provides a general method to detect all summands of a given subvariety without assuming minimality of their cohomology class. We hope that similar techniques may be useful for the question of Pareschi and Popa whether there are any other theta divisors with summands~\cite[conj.~19]{SchreiederCurveSummands} and for Debarre's minimal class conjecture~\cite{DebMinimal}. 

\medskip 

To find all decompositions of a given subvariety $Z\subset A$ as a sum $Z=X+Y$ of subvarieties $X,Y\subset A$, we 
use the theory of perverse sheaves and the decomposition theorem of Beilinson, Bernstein, Deligne and Gabber~\cite{BBD,DCM}. 
Let $\Perv(\bbC_A)$ be the category of perverse sheaves on the abelian variety and consider the intersection cohomology sheaves $\delta_X,\delta_Y,\delta_Z \in \Perv(\bbC_A)$ on $X,Y,Z\subset A$. Since these are simple perverse sheaves, the decomposition theorem for the addition map $a: X\times Y \twoheadrightarrow Z$ gives an embedding 
\[ \delta_Z \;\hookrightarrow \; \delta_X * \delta_Y \;=\; Ra_*(\delta_X \boxtimes \delta_Y) \]
as a direct summand in the derived category of constructible sheaves. Here $*$ denotes the additive convolution of perverse sheaves, which is controlled by the Tannakian formalism of~\cite{KraemerSemiabelian,KrWVanishing}. In particular, the above embedding induces an epimorphism
\[
 p: \quad G(\delta_X \oplus \delta_Y) \;\twoheadrightarrow\; G(\delta_Z)
\]
on the Tannakian groups of the respective perverse sheaves. The key point of our argument is that often this is an isogeny, so most information about the unknown summands $X,Y$ can be recovered from $Z$. Indeed, the results of~\cite{KraemerMicrolocalII} imply that up to isogeny these summands can be read off from the characteristic cycles of the perverse sheaves that correspond to the irreducible representations of the reductive group $G(\delta_Z)$, see theorem~\ref{thm:basicstep}.  For Jacobians of curves these perverse sheaves are the Brill-Noether sheaves introduced in~\cite{WeBN}, and in sections~\ref{sec:BNhyp} and~\ref{sec:BNnonhyp} we describe their characteristic cycles via Abel-Jacobi maps and highest weights for~$\Sp_{2g-2}(\bbC)$ and $\Sl_{2g-2}(\bbC)$. For intermediate Jacobians we instead use the representation theory of the exceptional group $E_6(\bbC)$, see section~\ref{sec:JT}.

\section{Reduction to Brill-Noether-Sheaves}

Let $A$ be a complex abelian variety and $Z\subset A$ an irreducible proper closed subvariety. In this section we describe a general method to determine all possible decompositions 
\[ 
 Z \;=\; X+Y 
\]
as a sum of irreducible geometrically nondegenerate subvarieties $X,Y\subset A$. For the notion of geometric nondegeneracy we refer to~\cite[p.~466]{RanSubvarieties}.

\subsection{Perverse sheaves and the decomposition theorem}
The first step is to reformulate such a decomposition in terms of the perverse intersection cohomology sheaves $\delta_X, \delta_Y, \delta_Z \in \Perv(\bbC_A)$:

\begin{lem} \label{lem:tensorproduct}
If $a: X\times Y \twoheadrightarrow Z$ denotes the addition map, we have an embedding as a direct summand
\[ \delta_Z \; \hookrightarrow \; Ra_*(\delta_X\boxtimes \delta_Y). \]
\end{lem}

{\em Proof.} The geometric nondegeneracy of the summands implies by~\cite[sect.~8.2]{DebarreComplexTori} that the addition morphism $a: X\times Y \twoheadrightarrow Z$ is generically finite and hence restricts to a finite \'etale cover over some smooth open dense subset $U\subseteq Z$.  Adjunction gives a direct summand~$\bbC_U \hookrightarrow a_*(\bbC_{X\times Y})|_U$, and shifting by $\dim Z = \dim X + \dim Y$ to make both sides perverse, we get the desired embedding over the smooth open dense subset~$U\subseteq Z$. The claim then follows by the decomposition theorem~\cite{BBD,DCM}. \qed 

\medskip 

The direct image in the above lemma can be rewritten as $Ra_*(\delta_X \boxtimes \delta_Y)=\delta_X*\delta_Y$ for the convolution product $*$ from~\cite{KraemerMicrolocal,KraemerSemiabelian,KrWVanishing}. To fix notations let us briefly recall from loc.~cit.~the Tannakian description of this convolution product.

\subsection{A Tannakian description} \label{subsec:convolution}

Generic vanishing on abelian varieties~\cite{KrWVanishing, SchnellHolonomic} says that for any given perverse sheaf $P\in \Perv(\bbC_A)$, its tensor product with a sufficiently general local system $L$ of rank one has the property
$ H^i(A, P\otimes L)=0$ for all $i\neq 0$.
We say that the perverse sheaf is {\em negligible} if it satisfies the following two equivalent conditions~\cite[th.~7.6]{SchnellHolonomic} \cite{ WeissauerVanishing}: \smallskip 
\begin{itemize} 
\item $H^0(A, P\otimes L)=0$ for a sufficiently general local system $L$ of rank one, \smallskip
\item each simple perverse subquotient in a Jordan-H\"older series for $P$ is invariant under translations by some positive-dimensional abelian subvariety. \smallskip
\end{itemize} 
Consider the abelian quotient category $\rmP(A)=\Perv(\bbC_A)/\rmS(A)$, where $\rmS(A)$ denotes the Serre subcategory of negligible perverse sheaves. It has been shown in~\cite{KrWVanishing,KraemerSemiabelian} that this quotient category is a rigid abelian tensor category with respect to the convolution product
\[ *: \quad \rmP(A)\times \rmP(A) \to \rmP(A),
\quad P_1*P_2 \;=\; Ra_*(P_1\boxtimes P_2)
\]
where $a: A\times A \to A, (x,y)\mapsto x+y$ is the addition map, and that for any~$P \in \rmP(A)$ the subcategory
\[ \langle P \rangle  \;=\; 
 \Biggl\{
 {\textnormal{smallest rigid abelian tensor}\atop\textnormal{subcategory containing $P$}}
 \Biggr\} \;\subset \; \rmP(A) 
\] 
is a neutral Tannakian category: It has a faithful exact $\bbC$-linear tensor functor~$\omega$ to the category $\Vect(\bbC)$ of finite dimensional complex vector spaces. Any such fiber functor induces an equivalence
\[
 \omega: \quad \langle P \rangle \; \stackrel{\sim}{\longrightarrow}\; \Rep(G)
\]
with the $\bbC$-linear abelian tensor category of complex algebraic representations of some affine algebraic group $G=G(P,\omega)$. Note that there is no canonical choice of the fiber functor; changing the fiber functor will replace the Tannakian group~$G$ by an isomorphic one. For an explicit choice one can proceed as in~\cite{GaL}, using that any~$Q\in \rmP(A)=\Perv(\bbC_A)/\rmS(A)$ has a unique representative $Q_\mathrm{int} \in \Perv(\bbC_A)$ without negligible subobjects or quotients: If~$L$ is a general local system of rank one, then
\[
  \langle P \rangle \; \longrightarrow\; \Vect(\bbC), \quad 
  Q \;\mapsto\; H^0(A, Q_\mathrm{int}\otimes L)
\]
is a fiber functor. Here the notion of a general local system depends on $P\in \rmP(A)$, but all arguments of this paper can be read in the Tannakian category generated by some fixed perverse sheaf, see section~\ref{subsec:summands-tannakian}. Any fiber functor on this ambient category gives by restriction a consistent choice of fiber functors on all its tensor subcategories. We fix such a consistent choice and simply write $G(P)=G(P, \omega)$ for any $P\in \rmP(A)$ that occurs in this paper. 

\begin{lem} \label{lem:abstractepi}
For any $Q\in \langle P\rangle$ we have a natural epimorphism $p: G(P)\twoheadrightarrow G(Q)$.
\end{lem}

{\em Proof.} By construction we have an embedding $\langle Q \rangle \hookrightarrow \langle P \rangle$ as a full tensor subcategory stable under subobjects and quotients, and we have chosen the fiber functors on these categories compatibly. So~\cite[cor.~2.9 and prop.~2.21]{DM} applies. \qed

\medskip 

\subsection{The ring of clean cycles} \label{subsec:cleancycles}

Before we come back to our concrete geometric problem, we need to recall a few facts about characteristic cycles. For a closed subvariety $Z\subset A$ we define its {\em conormal variety} $\Lambda_Z \subset T^*A$ as the Zariski closure of the conormal bundle to its smooth locus. This is a conic Lagrangian subvariety, and any conic Lagrangian subvariety $\Lambda \subset T^* A$ arises like this. By the {\em Gauss map} of $\Lambda$ we mean the projection \medskip 
\[ \gamma: \quad \Lambda \;\subset\; T^* A \;=\; A\times V \;\twoheadrightarrow \; V \;=\; H^0(A, \Omega^1_A), \quad (p,v) \;\mapsto \;v  \medskip \] 
to the fiber of the trivial cotangent bundle. This map is either generically finite or non-dominant. We denote its generic degree by $\deg(\Lambda)\in \bbN\cup \{0\}$ and say that~$\Lambda$ is {\em negligible} if $\deg(\Lambda)=0$. The degree extends additively to the group of conic Lagrangian cycles. Any cycle decomposes uniquely as a sum of a negligible and a clean part, where a cycle is called {\em clean} if it has no negligible components. The group
\[
\scrL(A) \;=\;  \{\textnormal{clean conic Lagrangian cycles on $T^*A$} \medskip
\}
\]
of clean cycles has a natural ring structure  where the product $\circ$ is induced by the correspondence
\[
\xymatrix@M=0.5em@C=6em{%
	T^* A \times T^*A \;=
	& A^2 \times V \ar[l]_-{\rho = \id \times diag} \ar[r]^-{\varpi=a\times \id}
	&  A\times V \;=\; T^* A
}
\]
coming from the addition map on the abelian variety~\cite{KraemerMicrolocalII}. By loc.~cit.~the map that sends a perverse sheaf to the clean part of its characteristic cycle is a ring homomorphism \medskip
\[
 \cc: \quad K^0(\rmP(A), \oplus, *) \;\longrightarrow\; (\scrL(A), +, \circ)
\] 
where the left hand side denotes the Grothendieck ring with respect to split exact sequences, with the ring structure given by convolution. We use the lowercase notation $\cc$ to emphasize that we only look at the clean part of the characteristic cycle, and similarly for the characteristic variety $\cha(-) = \Supp(\cc(-))$. 

\medskip 

For the proof of theorems~\ref{thm:JC} and~\ref{thm:JT} the ring structure is not enough, we need to control multilinear algebra constructions in the arising $\bbC$-linear abelian tensor categories. The Grothendieck ring $R$ of any such category is endowed with a family of maps $\lambda^n: R\to R$ induced by the alternating powers for $n\in \bbN_0$. These satisfy the relations
\[
 \lambda^0(a)=1, \quad \lambda^1(a)=a, \quad 
 \lambda^n(a+b)=\sum_{i=0}^n \lambda^i(a)\lambda^{n-i}(b)
\]
for all $a,b\in R$, and $\lambda^n(1)=0$ for all $n>1$. By~\cite[lemma~4.1]{Heinloth} we also have the identities
\begin{align*}
 \lambda^n(ab) &= P_n(\lambda^1(a), \dots, \lambda^n(a), \lambda^1(b), \dots, \lambda^n(b)), \\
 \lambda^m(\lambda^n(a)) &= P_{m,n}(\lambda^1(a), \dots, \lambda^{mn}(a)),
\end{align*}
where $P_{n}, P_{m,n}$ are universal polynomials with integer coefficients that describe the effect of composing alternating powers with tensor products resp.~with alternating powers. A commutative unital ring $R$ with maps $\lambda^n: R\to R$ satisfying all these properties is called a {\em $\lambda$-ring}. Rather than working directly with the maps $\lambda^n$ it is often more convenient to use the {\em Adams operations}. These are natural ring homomorphisms $\Psi^n: R\to R$ which can be expressed in terms of the $\lambda$-operations in the same way as the power sum polynomials are in terms of the elementary symmetric ones: We have
$\Psi^n=f_n(\lambda^1, \dots, \lambda^n)$,
where $f_n\in \bbZ[t_1, \dots, t_n]$ is defined by 
\[
 p_n \;=\; f_n(e_1, \dots, e_n) 
 \quad \textnormal{for} \quad 
 \begin{cases}
 \;p_n \;=\; \sum_{i=1}^\infty x_i^n \\
 \;e_n \;=\; \sum_{i_1<\cdots < i_n} x_{i_1} \cdots x_{i_n}.
 \end{cases}
\]
If $R$ has no $\bbZ$-torsion, then from the $\Psi^n$ for all $n\in \bbN$ we can recover the $\lambda^n$ since the elementary symmetric polynomials are rational polynomials in power sums.

\medskip 

\begin{ex} \label{ex:lambdaring}
The ring $\scrL(A)$ has no $\bbZ$-torsion. By~\cite{KraemerMicrolocalII} it is a $\lambda$-ring whose Adams operations 
\[
 \Psi^n \;=\; [n]_*: \quad \scrL(A) \;\longrightarrow\; \scrL(A) 
\]
are the pushforward of cycles under the map $[n]: A\times V \to A\times V, (x,v)\mapsto (nx,v)$, and 
\[
 \cc: \quad K^0(\rmP(A), \oplus, *) \;\longrightarrow\; (\scrL(A), +, \circ)
\]
is a homomorphism of $\lambda$-rings in the sense that we have $\cc(\Psi^n(P)) = [n]_* (\cc(P))$ for all~$ P \in \rmP(A)$. In what follows we will usually work in suitable finitely generated subrings of the above rings. For any cycle $\Lambda \in \scrL(A)$, we denote by 
$\langle \Lambda \rangle \subset \scrL(A)$
the smallest subring which contains this cycle and is stable under taking irreducible components of the support of cycles. By loc.~cit.~this subring is stable under the Adams operations, and restricting the homomorphism $\cc$ we get a homomorphism of $\lambda$-rings
\[ 
\cc: \quad K^0(\langle P \rangle, \oplus, *) \;\longrightarrow \; \langle \cc(P) \rangle \quad \textnormal{for any} \quad P\;\in\; \rmP(A).\]
\end{ex}

\pagebreak

\subsection{How to see summands on the Tannakian side} \label{subsec:summands-tannakian}

Let $Z=X+Y \subset A$ be a proper subvariety that decomposes as a sum of geometrically nondegenerate subvarieties of positive dimension. In the following we assume the summands $X,Y$ are fixed, though a priori not known explicitly. We want to reconstruct them from a given $P_Z\in \Perv(\bbC_A)$, which later will be taken to be the intersection cohomology sheaf on a translate of the curve or the Fano surface; see theorem~\ref{thm:basicstep}. In what follows we put
\[ \delta_{XY}=\delta_X \oplus \delta_Y \;\in\; \Perv(\bbC_A). \]
All that follows will take place inside the Tannakian category $\langle \delta_{XY}\oplus P_Z\rangle$. We fix a fiber functor $\omega$ on this ambient category and hence obtain a consistent choice of fiber functors for all the perverse sheaves to be considered below. Lemma~\ref{lem:tensorproduct} can then be reformulated as follows:

\begin{lem} \label{lem:epi}
We have a natural epimorphism $p: G(\delta_{XY})\twoheadrightarrow G(\delta_Z)$.
\end{lem}

{\em Proof.} Lemma~\ref{lem:tensorproduct} says that $\delta_Z \in \langle \delta_X * \delta_Y \rangle \subset \langle \delta_{XY} \rangle$, so lemma~\ref{lem:abstractepi} applies since we have chosen the fiber functors on our tensor categories in a consistent way. \qed

\medskip 

Note that $G(\delta_Z)\neq \{1\}$, indeed the stabilizer $\Stab(Z)=\{a\in A(\bbC) \mid Z+a = Z\}$ is finite since the sum of geometrically nondegenerate subvarieties is geometrically nondegenerate~\cite[cor.~8.11]{DebarreComplexTori}. 
Applying an isogeny to the abelian variety one reduces to the case where this stabilizer is trivial, and replacing all subvarieties by certain translates we can make their Tannakian groups semisimple~\cite[lemma~5.3.1]{KraemerMicrolocalII}. Let us now assume
\begin{align} \label{ass:stabilizer}  \tag{i}
 & \Stab(Z) \;=\; \{0\}, \\
 \label{ass:semisimple} \tag{ii}
 & G(\delta_{XY}) \; \textnormal{is a semisimple group}, \\
 \label{ass:gaussfinite} \tag{iii}
 & \textnormal{the Gauss map} \; \gamma: \; \bbP(\cha(\delta_Z)) \;\twoheadrightarrow\; \bbP V
 \; \textnormal{is a finite morphism}.
\end{align}
The last condition seems unreasonably strong but 
holds in our applications~\cite{KraemerThreefolds}.

\begin{prop} \label{prop:isogeny}
If~\eqref{ass:stabilizer} -- \eqref{ass:gaussfinite} hold, then $p:G(\delta_{XY})\twoheadrightarrow G(\delta_Z)$ is an isogeny.
\end{prop}

{\em Proof.} We must show that the kernel of the epimorphism $p$ is finite. On the level of Lie algebras this epimorphism splits as the projection to a direct summand. So we have
\[
 \Lie(G(\delta_{XY})) \;=\; \frakg_1 \oplus \frakg_2 
 \quad \textnormal{where} \quad 
 \begin{cases}
 \;\frakg_1 \;=\; \Lie(G(\delta_Z)), \\
 \;\frakg_2 \;=\; \Lie(\ker(p)),
 \end{cases}
\] 
and we want to show that the second summand in this decomposition vanishes. By assumption~\eqref{ass:stabilizer} we know
$\Stab(S) \subseteq \Stab(Z)=\{0\}$ for $S=X,Y$ and hence $\omega(\delta_S)$ remains an irreducible representation when restricted to the connected component of $G(\delta_{XY})$~\cite[cor.~1.6]{KraemerMicrolocal}. As a representation of the product $\Lie(G(\delta_{XY}))=\frakg_1 \times \frakg_2$ then
\begin{align}
  \label{eq:decomposition} \tag{iv}
 & \omega(\delta_S) \;\simeq\; V_{S,1}\boxtimes V_{S,2} \;\; \textnormal{with irreducible}\;\; V_{S,i}\in \Rep(\frakg_i),\, i=1,2,\\
 \label{eq:inclusion} \tag{v}
 & \omega(\delta_Z) \;\hookrightarrow\; (V_{X,1}\otimes V_{Y,1})\boxtimes \one \;\hookrightarrow\; (V_{X,1}\otimes V_{Y,1})\boxtimes (V_{X,2}\otimes V_{Y,2}). 
\end{align}
Now $\frakg_2$ is a semisimple Lie algebra by assumption~\eqref{ass:semisimple}, so in order to show that it is trivial we only need to see that it is abelian. Hence we will be done if $\dim V_{S,2} = 1$ for~$S=X,Y$ because $V_{X,2}\oplus V_{Y,2}\in \Rep(\frakg_2)$ is a faithful representation. To check that both dimensions are one we look at characteristic cycles. 

\medskip 

By~\cite[sect.~2.3]{KraemerMicrolocalII} the decomposition~\eqref{eq:decomposition} is mirrored on the level of characteristic cycles after some isogeny: There exists a natural number~$n\in \bbN$ such that
\begin{equation} \label{eq:cc-decomposition} \tag{vi}
 \cc(\delta_{nS}) \;=\; \Lambda_{S,1} \circ \Lambda_{S,2}
 \quad \textnormal{with effective cycles} 
 \quad 
 \begin{cases} 
 \Lambda_{S,1} \;\in\; \langle \cc(\delta_Z)\rangle \\
 \Lambda_{S,2} \;\in\; \langle \cc(\delta_{XY})\rangle
 \end{cases}
\end{equation}
of degree $\deg(\Lambda_{S,i})=\dim(V_{S,i})$. Similarly, the inclusion~\eqref{eq:inclusion} gives an inequality of effective cycles
\begin{equation} \label{eq:cc-inclusion} \tag{vii}
\cc(\delta_{nZ}) \;\leq\; \Lambda_{X,1}\circ \Lambda_{Y,1}.
\end{equation}
Here~$\circ$ denotes the product in the ring $\scrL(A)$ of section~\ref{subsec:cleancycles}, and for~$P\in \rmP(A)$ we consider  as in example~\ref{ex:lambdaring} the smallest subring $\langle \cc(P) \rangle \subset \scrL(A)$ containing $\cc(P)$ and stable under the passage to irreducible
components. Note that in~\eqref{eq:cc-decomposition} the conormal variety
$
 \Lambda_{nS} \;\leq\; \cc(\delta_{nS}) \;=\; \Lambda_{S,1}\circ \Lambda_{S,2}
$
enters with multiplicity one and is the unique component whose base is a subvariety of dimension $\dim(S)$ in~$A$. So in order to show $\deg(\Lambda_{S,2})=1$ we only need to see that~$\Lambda_{S,2}$ is a sum of conormal varieties to a finite number of points in $A$, since then there can be only one such point. In other words we want to show
\[
 \dim(\pi(\Supp(\Lambda_{S,2}))) \;=\; 0
 \quad \textnormal{for the projection} \quad \pi: \quad T^* A \;\to\; A.
\] 
To check this last claim we consider Chern-Mather classes. If $\Lambda \in \scrL(A)$ is any clean effective cycle and we put $d(\Lambda)=\dim(\pi(\Supp(\Lambda)))$, then by~\cite[lemma~3.1.2]{KraemerMicrolocalII} we know that
\[
 c_{M,i}(\Lambda) \;\in\; \CH_i(A) \;\textnormal{is} \;
 \begin{cases}
 \textnormal{zero for all $i>d(\Lambda)$}, \\
 \textnormal{a non-zero effective cycle for all $i \in \{0,1,\dots, d(\Lambda)\}$}.
 \end{cases}
\]
So we only need to show the vanishing $c_{M,1}(\Lambda_{S,2})=0$. Since $\Lambda_{S,1}\in \langle \cc(\delta_Z)\rangle$, the Gauss map
$
 \gamma: \bbP \Supp (\Lambda_{S,1}) \rightarrow \bbP V
$
is finite by assumption~\eqref{ass:gaussfinite}. By~\cite[lemma~3.3.1]{KraemerMicrolocalII} then
\begin{equation} \label{eq:pontryagin} \tag{viii}
 c_{M,d}(\Lambda_{nS}) \;=\; c_{M,d}(\Lambda_{S,1} \circ \Lambda_{S,2}) \;=\; \sum_{i=0}^d c_{M, d-i}(\Lambda_{S,1}) * c_{M,i}(\Lambda_{S,2})
\end{equation}
for $d<\dim A$ and the Pontryagin product $*$ on the Chow ring. The Pontryagin product of non-zero effective cycles is non-zero effective. Since for $d=\dim(S)+1$ the left hand side in~\eqref{eq:pontryagin} vanishes, all summands on the right hand side must vanish as well, so we get
\begin{align*} \label{eq:inequality} \tag{ix}
& d(\Lambda_{S,1}) + d(\Lambda_{S,2})
\;\leq\; \dim(S).
\end{align*}
Putting everything together we have
\begin{align*}
 \dim(Z) 
 & \leq d(\Lambda_{X,1}\circ \Lambda_{Y,1})
 &\textnormal{by inequality~\eqref{eq:cc-inclusion}} 
 \\
 & \leq d(\Lambda_{X,1}) +  d(\Lambda_{Y,1}) 
 &\textnormal{by definition of $\circ$}
 \\
 & \leq \dim(X) + \dim(Y) - d(\Lambda_{X,2}) - d(\Lambda_{Y,2}) 
 &\textnormal{by inequality~\eqref{eq:inequality}}
 \\
 & = \dim(Z) - d(\Lambda_{X,2}) - d(\Lambda_{Y,2}) 
 & \hspace*{-3em}\textnormal{since $\dim(Z)=\dim(X)+\dim(Y)$}
\end{align*}
which implies $d(\Lambda_{X,2})=d(\Lambda_{Y,2})=0$ as required. \qed


\subsection{Adams operations} 

Usually the above isogeny is not an isomorphism, but we can see representations of any finite cover of a connected semisimple group~$H$ as virtual representations in the ring $R(H)=K^0(\Rep(H), \oplus, \otimes)$:

\begin{lem} \label{lem:adams}
Let $p:G\twoheadrightarrow H$ be an isogeny of connected reductive groups. If $n\in \bbN$ is sufficiently divisible, then the $n$-th Adams operation on the representation ring of~$G$ factors as
\[ 
\xymatrix@M=0.5em{
 R(G) \ar[rr]^-{\Psi^n} \ar@{..>}[dr]_-{\exists !\, \Psi^n_p} && R(G) \\
 & R(H) \ar[ur]_-{\iota_p}
}
%
\]
where $\iota_p: R(H) \hookrightarrow R(G)$ is the embedding of $\lambda$-rings induced by the isogeny $p$.
\end{lem} 

{\em Proof.} If $T\subset G$ is a maximal torus with character group $\sfX = \CartierDual{T}$, the character map sending a representation to its weight space decomposition gives an embedding $R(G) \hookrightarrow \bbZ[\sfX]$. In these terms the Adams operation $\Psi^n$ is induced by the $n$-th power map of the torus~\cite[prop.~II.7.4]{BtD}. But any isogeny of tori factors over the $n$-th power map for $n$ dividing the degree of the isogeny. \qed

\medskip 

Applying this to the universal cover of the group $H=G(\delta_Z)$ when this cover is realized by a given $P_Z\in \Perv(\bbC_A)$, we recover the unknown summands $X, Y \subset A$ by the following

\begin{thm} \label{thm:basicstep}
Assuming~\eqref{ass:stabilizer} -- \eqref{ass:gaussfinite}, let $P_Z\in \Perv(\bbC_A)$ be given with~$G(P_Z)$ simply connected and 
\[ \delta_Z \;\in\; \langle P_Z \rangle, \]
then the associated tensor category contains simple perverse sheaves $\varepsilon_X, \varepsilon_Y \in \langle P_Z \rangle$ with
\[
 \cc(\delta_{nX}) = [n]_* \cc(\varepsilon_X)
 \quad \textnormal{\em and} \quad 
 \cc(\delta_{nY}) = [n]_* \cc(\varepsilon_Y)
 \quad \textnormal{\em for some $n\in \bbN$}.
\]
\end{thm}

{\em Proof.} As before we have $\Stab(S)=\{0\}$ and hence $\delta_{nS} = [n]_*\delta_{S}$ for $S=X,Y,Z$ and all $n\in \bbN$. So replacing all perverse sheaves by their direct image under an isogeny we may assume all Tannakian groups to be connected~\cite[cor.~1.6]{KraemerMicrolocal}.~The isogeny $p: G(\delta_{XY})\twoheadrightarrow G(\delta_Z)$ in proposition~\ref{prop:isogeny} is then dominated by the universal cover, which in our case is the epimorphism $q: G(P_Z) \twoheadrightarrow G(\delta_Z)$ coming from the inclusion $\delta_Z \in \langle P_Z\rangle$. Let $\varphi$ be the unique isogeny making the following diagram commute:
\[
\xymatrix@M=0.5em@C=4em{
 G(P_Z) \ar@{..>}[r]^-{\exists ! \, \varphi} \ar[d]_-q & G(\delta_{XY}) \ar[d]^-p \\
 G(\delta_Z) \ar[r]^\id &  G(\delta_Z)
}
\]
Let $\iota_\varphi: \langle \delta_{XY} \rangle \hookrightarrow \langle P_Z \rangle$ be the embedding of tensor categories induced by $\varphi$ on the level of representation categories. This $\iota_\varphi$ comes from abstract Lie theory and usually it is not the identity on perverse sheaves: It can happen that for $S=X,Y$ the simple perverse sheaves
\[ 
 \varepsilon_S \;=\; \iota_\varphi(\delta_S) \; \in \; \langle P_Z\rangle 
\] 
satisfy $\cc(\varepsilon_S)\neq \cc(\delta_S)$. But by construction we have the following commutative diagram
\[
\xymatrix@M=0.5em{
 \langle \delta_Z \rangle \ar@{=}[r] \ar@{^{(}->}[d]_-{\iota_p} 
 & \Rep(G(\delta_Z)) \ar[r]^-\id \ar[d]
 & \Rep(G(\delta_Z)) \ar@{=}[r] \ar[d]
 & \langle \delta_Z \rangle \ar@{^{(}->}[d]_-{\iota_q}
 \\
 \langle \delta_{XY} \rangle \ar@{=}[r]
 & \Rep(G(\delta_{XY})) \ar[r]^-{\iota_\varphi}
 & \Rep(G(P_Z)) \ar@{=}[r]
 & \langle P_Z \rangle
}
\]
where the outer vertical arrows are the inclusion functors on the level of perverse sheaves, and in contrast with the abstract embedding $\iota_\varphi$ these inclusion functors satisfy
\[
\cc(\iota_q(\psi)) \;=\; \cc(\psi) \;=\; \cc(\iota_p(\psi))
\quad \textnormal{for all} \quad \psi \;\in\; \langle \delta_Z \rangle.
\]
%
%
%
Altogether
\begin{align*}
 [n]_* \circ \cc \circ \iota_\varphi
 &\;=\; \Psi^n \circ \cc \circ  \iota_\varphi
 & \textnormal{since $\Psi^n = [n]_*$ on $\scrL(A)$ by example~\ref{ex:lambdaring}}
 \\
 &\;=\; \cc \circ \iota_\varphi \circ \Psi^n
 & \textnormal{by naturality of the Adams operations}
 \\
 &\;=\; \cc \circ \iota_\varphi \circ \iota_p \circ \Psi^n_p
 & \textnormal{by lemma~\ref{lem:adams} for the isogeny $p$}
 \\
 &\;=\; \cc \circ \iota_q \circ \Psi^n_p 
 & \textnormal{by functoriality for $q=p\circ \varphi$} 
 \\
 &\;=\; \cc \circ  \iota_p \circ \Psi^n_p 
 & \textnormal{since $\cc \circ \iota_q = \cc = \cc\circ \iota_p$ on $\langle \delta_Z \rangle$}
 \\
 &\;=\; \cc \circ \Psi^n
 & \textnormal{by lemma~\ref{lem:adams} for the isogeny $p$}
 \\
 &\;=\; \Psi^n \circ \cc
 & \textnormal{by naturality of the Adams operations}
 \\
 &\;=\; [n]_* \circ \cc
 & \textnormal{since $\Psi^n = [n]_*$ on $\scrL(A)$ by example~\ref{ex:lambdaring}}
\end{align*}
and it only remains to plug in the perverse sheaves $\delta_X, \delta_Y \in \langle \delta_{XY}\rangle$.
 \qed

\subsection{Highest weights}

The above reduces the classification of all geometrically nondegenerate summands of a given $Z\subset A$ to the study of the cycles~$\cc(\varepsilon)$ for the countably many simple perverse sheaves $\varepsilon \in \langle P_Z\rangle$. If~$P_Z$ is known, this is mainly a computational task using the highest weight theory of the connected reductive group 
\[ G \;=\; G(P_Z). \] 
In the next sections we will carry out the computations for Jacobians, but let us first fix some general notations. Let $T\subset G$ be a maximal torus with character group $\sfX=\CartierDual{T}$. Sending a representation to its character we obtain an isomorphism
\[ 
  \ch: \quad R(G) \;\stackrel{\sim}{\longrightarrow}\; \bbZ[\sfX]^W 
\] 
from the representation ring to the invariants of the Weyl group $W = N_G(T)/T$ in the character ring. Let us fix a system~$\Delta$ of simple positive roots and denote by~$\sfX^+\subset \sfX$ the corresponding set of dominant weights. Any Weyl group orbit contains a unique dominant weight, hence a natural additive basis of $\bbZ[\sfX]^W$ are the orbits
\[
 W e_\lambda \;=\; \sum_{w\in W} e_{w \lambda} \;\in\; \bbZ[\sfX]^W
 \quad \textnormal{for} \quad \lambda \;\in\; \sfX^+
\]
where $e_\mu \in \bbZ[\sfX]$ are the standard basis vectors of the group ring for $\mu \in \sfX$. On the other hand, the representation ring has a natural basis consisting of the irreducible representations $V_\lambda \in \Rep(G)$ of highest weight $\lambda$, where the latter again runs over all dominant weights. The transition between the two bases is given by the Weyl character formula which allows to compute the multiplicities $m_\lambda(\mu)\in \bbN_0$ of the weights in
\[
 \ch(V_\lambda) \;=\; \sum_{\mu\in \sfX^+} m_\lambda(\mu) \cdot We_\mu.
\]
We will only need to know which weights have positive multiplicity. This is easily answered as follows~\cite[cor.~3.2.12]{GW}, using the dominance order in which~$\mu \preceq \lambda$ iff~$\lambda - \mu$ is a sum of simple positive roots:

\begin{rem} \label{rem:weights}
For $\lambda, \mu \in \sfX^+$ one has $m_\lambda(\mu) > 0$ if and only if $\mu \preceq \lambda$.
\end{rem}

To any $\lambda \in \sfX^+$ we now attach via our chosen fiber functor $\omega: \langle P_Z \rangle \stackrel{\sim}{\longrightarrow}\Rep(G)$ a simple perverse sheaf
$\varepsilon_\lambda \in  \langle P_Z \rangle$ with $\omega(\varepsilon_\lambda) \simeq V_\lambda$. 
These form a complete set of representatives for the isomorphism classes of simple perverse sheaves in $\langle P_Z \rangle$ generalizing  the Brill-Noether sheaves from~\cite{WeBN}. To describe their characteristic cycles, let
\[
 \cc: \quad 
 \quad \bbZ[\sfX]^W \;\simeq\; R(G) \;\simeq\; \langle P_Z \rangle \;\longrightarrow\; \langle \cc(P_Z) \rangle \;\subset\; \scrL(A) 
\]
denote the composition of the characteristic cycle with the inverse of $\ch \circ \omega$. Then we have
\[
 \cc(\varepsilon_\lambda) \;=\;
 \sum_{\mu \in \sfX^+} m_\lambda(\mu) \cdot \cc(\mu)
 \quad \textnormal{where} \quad 
 \cc(\mu) \;=\; \cc(W e_\mu),
\]
hence we know the possible supports in theorem~\ref{thm:basicstep} if we can control $\cc(\mu)$ for all dominant weights~$\mu$. The cases to be considered below are particularly simple because here the monodromy action on a general fiber of the Gauss map coincides with the Weyl group action on weights~\cite[th.~9]{KraemerThreefolds}~\cite[th.~2.1]{KraemerMicrolocal}:

\begin{rem} \label{rem:weyl-monodromy}
If $P_Z$ is the perverse intersection cohomology sheaf of a curve in its Jacobian or of the Fano surface of a smooth cubic threefold in the intermediate Jacobian, then 
\[ \cc: \quad \bbZ[\sfX]^W \;\stackrel{\sim}{\longrightarrow} \; \langle \cc(P_Z) \rangle \]
is an isomorphism and all the cycles $\cc(\mu)$ are reduced and irreducible. \end{rem}
 
Recall that every closed conic Lagrangian subvariety of the cotangent bundle is the conormal variety to its base. An explicit description of $\cc(\mu) \subset T^* A$ as a conormal variety will be given in lemmas~\ref{lem:easyorbithyp} and~\ref{lem:easyorbit}.

\section{Brill-Noether sheaves on hyperelliptic Jacobians} \label{sec:BNhyp}

Let $A=JC$ be the Jacobian variety of a smooth curve of genus $g>1$. In this section we assume the curve to be hyperelliptic and embed it as a symmetric subvariety
$
 C = -C \subset A
$
via some translate of the Abel-Jacobi map. Then $G(\delta_C)$ is the full symplectic group of rank $n=g-1$~\cite[th.~6.1]{KrWSmall}~\cite{WeBN}, and we have a fiber functor
\[ 
  \omega: \quad \langle \delta_C \rangle  \; \stackrel{\sim}{\longrightarrow} \; \Rep(\Sp_{2n}(\bbC)) 
\] 
where $\omega(\delta_C)=\bbC^{2n}$ is the natural representation of the group $\Sp_{2n}(\bbC)$.

\subsection{Weights for the symplectic group}

To fix notations we consider the split form
\[
 \Sp_{2n}(\bbC) 
 \;=\;
 \bigl\{
 M \in \Gl_{2n}(\bbC) \mid
 M^t \, \Omega\,  M = \Omega
 \bigr\} 
\]
where
\[
 \Omega \;=\; 
 \Biggl(
 \begin{matrix}
 \;\;0 & \mathbb{I} \\
 -\mathbb{I} & 0
 \end{matrix} \;
 \Biggr)
 \quad \textnormal{for the $n\times n$ matrix} \quad 
 \mathbb{I} \;=\;
 \Biggl(
 \begin{smallmatrix}
  0 & & 1 \\
    & \updots & \\
  1 & & 0
 \end{smallmatrix}
 \Biggr) 
\]
so that the set of diagonal matrices
$
 T = 
 \bigl\{
 \, \mathrm{diag}(t_1, \dots, t_{2n}) 
  \mid t_{n+i} = t_i^{-1} \, \textnormal{for all} \, i \,
 \bigr\}
$
is a maximal torus. Sending such a diagonal matrix to its $i$-th entry defines a character $e_i\in \sfX=\CartierDual{T}$ for $1\leq i\leq n$. These form a basis of the lattice $\sfX$ in which we choose
\begin{align*}
 \Delta
 &\;=\; \{ e_i - e_{i+1} \mid i = 1, \dots, n-1\} \cup \{ 2e_n\}, \\
 \sfX^+ 
 &\;=\; \{ \lambda \in \bbZ^n \mid \lambda_1 \geq \cdots \geq \lambda_n \geq 0 \} 
 \;\subset\; \sfX \;=\; \bbZ^n.
\end{align*}
Another basis of the weight lattice are the fundamental weights $\varpi_d = e_1 + \cdots + e_d$ for~$1\leq d\leq n$, which are important since they generate the semigroup of dominant weights. The fundamental representations $V_{\varpi_d} \in \Rep(\Sp_{2n}(\bbC))$ appear in the alternating powers
\[
 \Alt^d(\bbC^{2n})
 \;\simeq\; \bigoplus_{i=0}^{\lfloor d/2 \rfloor} V_{\varpi_{d-2i}}
 \quad \textnormal{for} \quad
 d \;=\; 1,\dots, n. 
\]
On perverse sheaves this is reflected by the geometry of the Abel-Jacobi map:

\begin{ex}
For $d=1,\dots, n$ consider the sum map $a_d: C_d = C^d/\frakS_d \longrightarrow A$ for our chosen translate of the curve. In the hyperelliptic case a look at fiber dimensions shows that
\[
 Ra_{d*}(\delta_{C_d}) \;\simeq\; \bigoplus_{i=0}^{\lfloor d/2 \rfloor} \; \delta_{W_{d-2i}}
\]
and the Brill-Noether sheaves corresponding to the fundamental weights are the perverse intersection cohomology sheaves
$
 \varepsilon_{\varpi_d} \simeq \delta_{W_d}
$~\cite[lemma~20]{WeBN}. 
\end{ex}

Note that the supports in this example are precisely the summands from the main theorem. In order to generalize this to arbitrary dominant weights $\lambda \in\sfX^+ \subset \bbZ^{n}$ we denote by
\begin{align*}
 d(\lambda) &\;=\; \lambda_1 + \cdots + \lambda_n, \\
 \ell(\lambda) &\;=\; \max\{1\leq i\leq n\mid  \lambda_i > 0\}\cup \{0\}
\end{align*}
the degree and length as a partition. We will apply remark~\ref{rem:weights} as follows:

\begin{lem} \label{lem:muhyp}
With notations as above, for any $\lambda \in\sfX^+$ there exists $\mu \in\sfX^+$ such that 
\[
 \mu \;\preceq \; \lambda
 \quad \textnormal{\em and} \quad 
 \ell(\mu) \;=\; \min\{ d(\lambda), n\}.
\]
\end{lem}

{\em Proof.} Clearly $\ell(\lambda) \leq \min\{ d(\lambda), n\}$. If equality holds, then taking $\mu = \lambda$ we are done, so let us assume that the inequality is strict. In this case we can define $i, k$ by
$
 1\leq i = \max\{ j \mid \lambda_j \geq 2 \} < k = \ell(\lambda) + 1 \leq n. 
$
Then
$
 \nu = \lambda - (e_i-e_k) \preceq \lambda
$
is a dominant weight with $d(\nu)=d(\lambda)$ but~$\ell(\nu)=\ell(\lambda) + 1$, so we are done by ascending induction on the length.  \qed

\medskip 

\subsection{Weyl group orbits and conormal varieties}

By remark~\ref{rem:weyl-monodromy} the irreducible components of characteristic cycles of Brill-Noether sheaves are precisely the conic Lagrangian subvarieties $\cc(\mu) \subset T^* A$. To describe them, we consider for $\mu\in \sfX^+$ the image
\[
 W(\mu) \;=\;
 \Im 
 \Bigl(
 a_\mu: 
 C^{n} \rightarrow A,
 \;\; p\mapsto \mu_1 p_1 + \cdots + \mu_n p_n
 \Bigr)
 \;\subseteq\; A.
\]
This image with the reduced subscheme structure is geometrically nondegenerate of dimension $\dim(W(\mu)) = \ell(\mu)$ by~\cite[cor.~8.11]{DebarreComplexTori}, and we claim that $\cc(\mu)$ is its conormal variety:

\begin{lem} \label{lem:easyorbithyp}
We have $\cc(\mu) = \Lambda_{W(\mu)}$. 
\end{lem}

{\em Proof.} Since $\omega(\delta_C)$ is the natural representation of the symplectic group $\Sp_{2n}(\bbC)$, its highest weight is the fundamental weight $\varpi_1 = e_1 = (1,0,\dots, 0)$, so the Weyl group orbit
\[ 
 \cc(\varpi_1) \;=\; \cc(\delta_C) \;=\; \Lambda_C
\]
is the conormal variety to the curve. Via the natural identification $V=H^0(C, \omega_C)$ we regard $\bbP V = |\omega_C|$ as the linear series of canonical divisors on the curve, in which case
\[
 \bbP \Lambda_C \;=\; \{ (p, D) \in C \times |\omega_C| \mid p \in \Supp(D) \}
\]
becomes the incidence variety for the divisors in the canonical linear series.
Now let $u\in V\setminus \{0\}$ be a non-zero global holomorphic differential form. Any canonical divisor is stable under the hyperelliptic involution $\iota\in \Aut(C)$, hence we may write it as
\[
 div(u) \;=\; p_1 + \cdots + p_n + \iota(p_1) + \cdots + \iota(p_n)
 \quad \textnormal{with} \quad p_1, \dots, p_n \;\in\; C \;\subset\; A,
\]
and then the group
\[
 \Gamma \;=\; \langle a\in A(\bbC) \mid (a, u) \in \Lambda_C \rangle
\]
is generated by $p_1, \dots, p_n$. On the other hand, for very general $u$ we know from~\cite{KraemerMicrolocal} that $\Gamma$ is a free abelian group of rank $n$, and hence in loc.~cit.~we can choose the isomorphism
\[
 p: \quad \sfX \;\stackrel{\sim}{\longrightarrow} \; \Gamma
 \quad \textnormal{such that} \quad e_i \;\mapsto \; p_i
\] 
where $e_i \in X$ are the standard basis weights from the previous section. Varying $u$ one deduces that for any $\mu \in\sfX^+$ the corresponding conic Lagrangian subvariety is given by
\[
 \bbP \cc(\mu) \;=\;
 \mathrm{Im}\bigl((a_\mu, q): C^n \rightarrow A\times \bbP V \bigr) 
\]
where the image on the right hand side is considered as a reduced closed subscheme and
$
  q:  C^n \twoheadrightarrow \bbP V = |\omega_C|,  
 D \mapsto D + \iota(D) 
$
is the map parametrizing canonical divisors as above. So for the projection from the characteristic cycle to the abelian variety it follows that $W(\mu) = \mathrm{Im}(\bbP \cc(\mu) \to A)$, hence $\cc(\mu)=\Lambda_{W(\mu)}$ as any irreducible Lagrangian cone in the cotangent bundle is the conormal variety to its base. 
\qed

\subsection{Characteristic cycles of Brill-Noether sheaves}

For hyperelliptic curves the above gives a complete description of the characteristic cycle of Brill-Noether sheaves. In particular we get

\begin{prop} \label{prop:easyBNhyp}
If $d=d(\lambda)\leq n$, then $\Supp(\varepsilon_\lambda)=p(\cha(\varepsilon_\lambda))=W_{d}$.
\end{prop}

{\em Proof.}
Any dominant weight can be written as a nonnegative linear combination of the fundamental weights $\varpi_i = e_1 + \cdots + e_i$ with $1\leq i\leq n$, hence $\lambda$ can be written as 
\[
 \lambda \;=\; \varpi_{i_1} + \cdots + \varpi_{i_m}
 \quad \textnormal{with} \quad 
 1\leq i_1 \leq \cdots \leq i_m \leq n.
\]
So by representation theory the irreducible representation of highest weight $\lambda$ enters in 
$\Alt^{i_1}(\bbC^{2n}) \otimes \cdots \otimes \Alt^{i_m}(\bbC^{2n})$~\cite[th.~5.5.21]{GW},
hence in the $d$-th tensor power of the natural representation. Geometrically this provides an embedding $\varepsilon_\lambda \hookrightarrow (\delta_C)^{*d}$ so that
$\Supp(\varepsilon_\lambda) \subseteq W_d$.
It remains to show that the conormal variety to the right hand side is an irreducible component of $\cc(\varepsilon_\lambda)$. By assumption $d=d(\lambda)\leq n$, so lemma~\ref{lem:muhyp} says that
\[
 \mu \;=\; (1^d, 0^{n-d}) \;\preceq \; \lambda.
\]
Since $W(\mu) = W_d$ by definition, we are done by remark~\ref{rem:weights} and lemma~\ref{lem:easyorbithyp}. \qed

\medskip 

Note that in the situation of the above corollary one has $\dim(W_\lambda) < \dim(W_d)$ unless $\lambda = (1^d, 0^{n-d})$. So in most cases the Weyl group orbit of the highest weight corresponds to a component of $\cc(\lambda)$ which does {\em not} dominate $\Supp(\varepsilon_\lambda)$.

\begin{thm} \label{thm:BNhyp}
For $\lambda \in\sfX^+$ the dimension of the subvariety $S = p(\cha(\varepsilon_\lambda))\subset A$ is given by
\[
 \dim S \;=\; \min \{d(\lambda), g-1\}. \medskip
\]
Hence if $S$ has dimension $d<g-1$, then it follows that $S=W_d$.
\end{thm}

{\em Proof.} By proposition~\ref{prop:easyBNhyp} we may assume that $d(\lambda) > n$. Combining remark~\ref{rem:weights} with lemma~\ref{lem:easyorbithyp} it will suffice to show that in this case there exists a dominant integral weight $\mu \preceq \lambda$ with $\ell(\mu)=n$. But this follows from lemma~\ref{lem:muhyp}. \qed

\subsection{Proof of the main theorem for hyperelliptic Jacobians} \label{sec:proofhyp}

If $\Theta = X+Y$ is a decomposition as a sum of two geometrically nondegenerate subvarieties of positive dimension, then theorem~\ref{thm:basicstep} says that there exists $n\in \bbN$ and $\lambda, \mu \in\sfX^+$ such that
\[
 \cc(\delta_{nX}) \;=\; [n]_* \cc(\varepsilon_\lambda)
 \quad \textnormal{and} \quad 
 \cc(\delta_{nY}) \;=\; [n]_* \cc(\varepsilon_\mu).
\]
Projecting to the abelian variety we obtain from the last statement in theorem~\ref{thm:BNhyp} that
\begin{align*}
 nX &\;=\; nW_d,
 \quad  d \;=\; \dim(X), \\
 nY &\;=\; nW_e,
 \quad \, e \;=\; \dim(Y).
\end{align*}
Taking preimages of these identities under the isogeny $[n]: A\longrightarrow A$ we deduce that \medskip 
\[
 W_d \;\subseteq\;
 \bigcup_{x\in A[n]} \bigl(X - x\bigr) 
 \quad \textnormal{and} \quad 
 W_e \;\subseteq\;
  \bigcup_{y\in A[n]} \bigl(Y - y\bigr). \medskip
\]
By irreducibility then $X=W_d+x$ and $Y=W_e+y$ for certain $x,y\in A[n]$. \qed

\medskip

\section{Brill-Noether sheaves on nonhyperelliptic Jacobians} \label{sec:BNnonhyp}

Now let $A=JC$ be the Jacobian variety of a smooth curve of genus $g>1$ that is not hyperelliptic. Here the Abel-Jacobi curve cannot be taken symmetric but we choose a translate so that for $n=g-1$ the divisor $W_n=-W_n\subset A$ becomes symmetric. Then $G(\delta_C)$ is a special linear group by~\cite[th.~6.1]{KrWSmall}~\cite{WeBN} and we have a fiber functor
\[
 \omega: 
 \quad \langle \delta_C \rangle 
 \;\stackrel{\sim}{\longrightarrow}\;
  \Rep(\Sl_{2n}(\bbC))
\]
where $\omega(\delta_C)=\bbC^{2n}$ is the natural representation of the group $\Sl_{2n}(\bbC)$.

\subsection{Weights for the special linear group}

Let $T\subset \Sl_{2n}(\bbC)$ be the maximal torus of diagonal matrices in the special linear group. The characters $e_i\in \sfX$ sending a diagonal matrix to its $i$-th entry generate the character lattice, with the relation~$e_1 + \cdots + e_{2n} = 0$. So in this case $\sfX = \bbZ^{2n}/\langle \det\rangle$ for $\det = (1,1,\dots, 1)$ and we choose
\begin{align*}
 \Delta
 &\;=\; \{ e_i - e_j \mid 1\leq i < j \leq 2n\}, \\
\sfX^+ 
 &\;=\; \{ \lambda \!\!\!\! \mod \langle \det \rangle \mid \lambda_1 \geq \cdots \geq \lambda_{2n} \geq 0 \}. 
\end{align*}
The fundamental weights 
$\varpi_d = e_1 + \cdots + e_d$ for~$d = 1,\dots, 2n-1$
are a basis of the character lattice and generate the semigroup of dominant weights, but now we have
\[
 V_{\varpi_d} \;\simeq\; \Alt^d(\bbC^{2n}) \;\simeq\; \Alt^{2n-d}(\bbC^{2n})^* \;\simeq\; \Hom(V_{\varpi_{2n-d}}, \one).
\]
Again this has a geometric interpretation:

\begin{ex}
In the nonhyperelliptic case the fundamental weights $\varpi_d \in\sfX^+$ correspond to 
\[ \varepsilon_{\varpi_d} \;=\;
\begin{cases}
 \delta_{W_d} & \textnormal{if $0\leq d\leq n$}, \\
 \delta_{-W_{n - d}} & \textnormal{if $g\leq d< 2n$},
\end{cases}
\]
as one may see via the Riemann-Roch theorem and Serre duality~\cite{WeBN}.
\end{ex}

With this example in mind we fix a system of representatives modulo $\det$ as follows. For $\alpha, \beta\in \bbN_0^n$ put
$(\alpha \mid -\beta) = ( \alpha_1, \dots, \alpha_n, -\beta_1, \dots, -\beta_n)$. Any $\lambda \in \bbZ^{2n}$ satisfies
\[
 \lambda \; \equiv \; (\lambda^+ \mid -\lambda^-)
 \!\!\!\mod \langle \det \rangle
 \quad
 \textnormal{for some} \quad 
 \lambda^\pm \;\in\; \bbN_0^n
\]
and the latter are determined uniquely by imposing that at least one entry of $\lambda^-$ vanishes. We put
\begin{align*} 
 d^\pm(\lambda) &= \lambda^\pm_1 + \cdots + \lambda^\pm_n,  
 & d(\lambda) = d^+(\lambda) + d^-(\lambda), \\
 \ell^\pm(\lambda) &= \bigl| \{i \mid \lambda_i^\pm \neq 0\} \bigr|,  
 & \ell(\lambda) = \ell^+(\lambda) + \ell^-(\lambda),
 \end{align*} 
so that in the previous example
$\dim \Supp(\varepsilon_{\varpi_d}) = d(\varpi_d)$ 
without case distinctions. 

\begin{lem} \label{lem:munonhyp}
For any $\lambda \in\sfX^+$ there exists a dominant weight $\mu \in\sfX^+$ with $\mu \preceq \lambda$ such that
\[
 \ell(\mu) \;=\; \min\{ d(\lambda), n\} 
 \quad \textnormal{\em or} \quad 
 \ell(\mu) \;=\; d(\mu) \;=\; n-1.
\]
\end{lem}


{\em Proof.} We may clearly assume that $\ell(\lambda) \neq \min\{d(\lambda), n\}$. Now there are two cases: If $\ell(\lambda) < \min\{d(\lambda), n\}$,  then as in lemma~\ref{lem:muhyp} we successively subtract simple positive roots 
\[ (e_i-e_k\mid 0) \quad \textnormal{or} \quad (0\mid e_i - e_k)
\quad \textnormal{with} \quad 1\;\leq\; i \;<\;k \;\leq \; n 
\] 
to find a dominant integral weight $\mu \preceq \lambda$ with $\ell(\mu) = \min\{ d(\lambda), n\}$. It remains to treat the case $\ell(\lambda) > \min\{d(\lambda), n\}$. In this second case $n < \ell(\lambda) \leq d(\lambda)$, and we successively subtract simple positive roots
\[ (e_i \mid -e_k) \quad \textnormal{with} \quad 1 \;\leq \; i, k \;\leq \; n \]
until we arrive at a dominant integral weight $\nu \preceq \lambda$ satisfying one of the following conditions:
\begin{itemize}
\item $d(\nu) \in \{n, n-1\}$, 
\item $d(\nu) > n$ but $\nu^- = 0$,
\item $d(\nu) > n$ but $\nu^+ = 0$.
\end{itemize}
Under each of these three conditions an argument as in the first case allows to find a dominant integral weight $\mu \preceq \nu$ with $\ell(\mu) = n$ or $\ell(\mu)=d(\mu)=n-1$. \qed

\subsection{Weyl group orbits and conormal varieties}

For $\mu \in \sfX^+$ consider now the image
\[
 W(\mu) \;=\;
 \Im \Bigl(
 a_\mu: \;
 C^{n}\times C^{n} \longrightarrow A, \;
 (p, q) \mapsto  \sum_{i=1}^n \,  \bigl( \mu^+_i\cdot  p_i - \mu^-_i\cdot q_i \bigr)
 \Bigr)
 \;\subseteq\; A,
\]
viewed as a reduced subscheme of $A$. This image is geometrically nondegenerate of dimension
$\dim(W(\mu)) = \min\{ \ell(\mu), g\}$ 
by~\cite[cor.~8.11]{DebarreComplexTori}, and $\cc(\mu)$ is its conormal variety whenever possible: 

\begin{lem} \label{lem:easyorbit}
If $\ell(\mu)<g$, then $\cc(\mu) = \Lambda_{W(\mu)}$.
\end{lem}

{\em Proof.} Since $\omega(\delta_C) = \bbC^{2n}$ is the natural representation of the group $\Sl_{2n}(\bbC)$, its highest weight is the fundamental weight $e_1 = (1,0,\dots, 0)$, so for this weight we know that
\[
 \cc(e_1) \;=\; \cc(\delta_C) \;=\; \Lambda_C
\]
is the conormal variety of the curve. Then as in the proof of lemma~\ref{lem:easyorbithyp}  it follows that for any~$\mu \in\sfX^+$ the fiber of $\cc(\mu) \subset A\times V$ over a very general $u\in V=H^0(C, \Omega_C^1)$ is 
\[
 \{ a\in A(\bbC) \mid (a, u) \in \Supp(\cc(\mu)) \} \;=\; \{ a_\mu(D) \mid D\in q^{-1}(div(u)) \},
\]
where $q$ is the quotient map to the symmetric power of the curve in the following diagram:
\[
\xymatrix@M=0.8em{
\; D \;\in\; q^{-1}|\omega_C| \ar@{^{(}->}[r] \ar@<5ex>[d] & C^{2g-2} \ar[r]^-{a_\mu} \ar[d]^-q & A\\
div(u)  \;\in\; |\omega_C| \ar@{^{(}->}[r] & C^{2g-2}/\frakS_{2g-2} &
}
\]
As $u\in V$ varies, the divisor $div(u)$ runs through the canonical linear series $|\omega_C|$ and hence any $g-1$ points of it can be moved independently from each other. So the composite map
\[ 
 q^{-1}|\omega_C| \;\hookrightarrow\; C^{2g-2} \twoheadrightarrow  C^{|I|}, \quad 
 p \;\mapsto\; (p_{i_1}, p_{i_2}, \dots)
\]
is dominant for any index set $I \subset \{1,2,\dots, 2g-2\}$ of cardinality $|I|<g$. It follows that
\[
 a_\mu(q^{-1}|\omega_C|) \;=\; a_\mu(C^{2g-2}) \;=\; W(\mu)
\]
for $\ell(\mu)<g$, and hence we are done. \qed

\medskip 

For $\ell(\mu)\geq g$ the above proof still shows that the subvariety $\cc(\mu) \subset T^* A$ is the conormal variety to the image $a_\mu(q^{-1}|\omega_C|) \subset A$, but it seems unclear how to describe this image in general. For instance, is its dimension always $g-1$?

\subsection{Characteristic cycles of Brill-Noether sheaves} 

Although lemma~\ref{lem:easyorbit} only applies to a certain range of weights, it suffices to prove the main theorem in the nonhyperelliptic case. As before we begin with small degrees:

\begin{prop} \label{prop:easyBN}
If $d(\lambda) < g$, then 
$\Supp(\varepsilon_\lambda) = W_{d^+(\lambda)} - W_{d^-(\lambda)}$
and in this case the support is a geometrically nondegenerate subvariety of dimension $d(\lambda)$.
\end{prop}

{\em Proof.} Any dominant weight can be written as a nonnegative linear combination of the fundamental weights $\varpi_i = e_1 + \cdots + e_i$ with $1\leq i\leq 2n-1$, hence $\lambda$ can be written as a sum
\[
 \lambda \;=\; \varpi_{\mu_1} + \cdots + \varpi_{\mu_m}
\]
with $\mu_1, \dots, \mu_m \in \bbN_0$ not necessarily distinct. With these notations we know from representation theory that the irreducible representation of highest weight $\lambda$ is a summand in 
\[
 \bigotimes_{i=1}^m \Alt^{\mu_i}(U) \;\simeq \; \bigotimes_{\mu_i \leq n} \Alt^{\mu_i}(U) \otimes \bigotimes_{\mu_i > n} \Alt^{2n-\mu_i}(U^*),
\]
where $U=\bbC^{2n}$ denotes the natural representation of $\Sl_{2n}(\bbC)$~\cite[th.~5.5.21]{GW}. In our previous notations
\begin{itemize}
\item the $\mu_i \leq n$ are the parts of the transpose partition $(\lambda^+)^t$
\item the $2n-\mu_i$ with $\mu_i > n$ are the parts of the transpose partition $(\lambda^-)^t$
\end{itemize} 
and so the respective irreducible representation enters in $U^{\otimes d^+(\lambda)}\otimes (U^*)^{\otimes d^-(\lambda)}$ where $d^\pm(\lambda)$ denotes the degree of $\lambda^\pm$ as above. Translating this back to geometry we get
\[
 \varepsilon_\lambda \;\hookrightarrow\; 
 \underbrace{\delta_C * \cdots * \delta_C}_{d^+(\lambda)} * \underbrace{\delta_{-C}*\cdots * \delta_{-C}}_{d^-(\lambda)}
\]
and therefore
\[
 \Supp(\varepsilon_\lambda) \;\subseteq\; W_{d^+(\lambda)} - W_{d^-(\lambda)}.
\] 
It only remains to show that the conormal variety to the right hand side enters as an irreducible component of the clean cycle $\cc(\varepsilon_\lambda)$. Since by assumption $d(\lambda)\leq n$, lemma~\ref{lem:munonhyp} shows that
\[
 \mu \;=\; (\underbrace{1, \dots, 1}_{d^+(\lambda)}, 0, \dots, 0, \underbrace{-1, \dots, -1}_{d^-(\lambda)}) \;\preceq \; \lambda.
\]
Since $W_\mu = W_{d^+(\lambda)}-W_{d^-(\lambda)}$, we are done by remark~\ref{rem:weights} and lemma~\ref{lem:easyorbit}. \qed

\begin{thm} \label{thm:BN}
For any $\lambda \in\sfX^+$ the dimension of the image $S=p(\cha(\varepsilon_\lambda)) \subset A$ satisfies
%
\[
 \dim S \;\geq \; \min \{d(\lambda), g-2\}. \medskip
\]
Hence if $d=\dim S<g-2$, then $S=W_e-W_{d-e}$ for some $e\in \{0,1,\dots, d\}$. 
\end{thm}

{\em Proof.} By proposition~\ref{prop:easyBN} we may assume that $d(\lambda) > n$. Combining remark~\ref{rem:weights} with lemma~\ref{lem:easyorbit} it will suffice to show that in this case there exists a dominant integral weight $\mu \preceq \lambda$ with $\ell(\mu)\in \{n,n-1\}$. But this follows from lemma~\ref{lem:munonhyp}. \qed

\subsection{Proof of the main theorem for nonhyperelliptic Jacobians}

Suppose as before that $\Theta = X+Y$ is a sum of geometrically nondegenerate subvarieties of positive dimension. Again by theorem~\ref{thm:basicstep} there exist $n\in \bbN$ and $\lambda, \mu \in\sfX^+$ such that
\[
 \cc(\delta_{nX}) = [n]_* \cc(\varepsilon_\lambda) 
 \quad \textnormal{and} \quad 
 \cc(\delta_{nY}) = [n]_* \cc(\varepsilon_\mu).
\]
If one of the two summands is a curve, we are done by the result of Schreieder~\cite{SchreiederCurveSummands}, so we assume 
\[ 
  d \;=\; \dim X \;<\; g-2
  \quad \textnormal{and} \quad 
  e \;=\; \dim Y \;<\; g-2. 
\]
We are then in the range of the last statement in theorem~\ref{thm:BN}. As in section~\ref{sec:proofhyp} we deduce
\begin{align*}
 X &\;=\; W_{a}-W_{d-a}+x \quad
 \textnormal{with} \quad x\in A[n], \; a\in \{0,1, \dots, d\}, \\
 Y &\;=\; W_{b}-W_{e -b}\,+\,y \quad  
 \textnormal{with} \quad y\in A[n], \; b\in \{0,1,\dots, e\}.
\end{align*}
It remains to show that either $a=b=0$ or $d-a=e-b=0$.
Putting $z=x+y$ we know from the above decompositions that the theta divisor can be written as a difference
\[
 W_{g-1} \;=\; W_{c} - W_{g-1-c} + z
 \quad \textnormal{with} \quad 
 c \;=\; a+b \;\in\; \{0,1,\dots, g-1\},
\]
and we need to show that
\[
 c\;\in\; \{0, g-1\}.
\]
One way to see this is to use the bijection between components of characteristic cycles and Weyl group orbits of dominant weights in remark~\ref{rem:weyl-monodromy}, and to observe that for $c\notin\{0, g-1\}$ the weights $\varpi_{g-1}$ and $\varpi_c + \varpi_{g-1+c}$ are not in the same Weyl group orbit. A more direct geometric argument would be to observe that the identity $W_{g-1} = W_c - W_{g-1-c} + z$ implies  by~\cite[lemma 1c]{MartensSpecial} or by~\cite[lemma~17]{SchreiederCurveSummands} that
\[ W_c \;=\; z-W_c, \] 
so $c\in \{0,g-1\}$ by Martens' theorem for nonhyperelliptic curves. \qed 

\medskip 

\section{Intermediate Jacobians of cubic threefolds} \label{sec:JT}

Let us finally consider the case of intermediate Jacobian $A=JT$ where $T\subset \bbP^4$ is a smooth cubic threefold. Let $S\subset A$ be the Fano surface of lines on the threefold which by~\cite{KraemerThreefolds} we can embed in the intermediate Jacobian via some translate of the Albanese map so that $G(\delta_S)$ is the simply connected group $E_6(\bbC)$. We have a fiber functor
\[
 \omega: \quad \langle \delta_S \rangle \;\stackrel{\sim}{\longrightarrow} \; \Rep(E_6(\bbC))
\]
where $\omega(\delta_S)$ is one of the two irreducible representations of dimension $27$.

\subsection{Weights for the group $E_6(\bbC)$}

We follow the conventions of Bourbaki~\cite{BourbakiLie} and label the simple positive roots by $\alpha_1, \dots, \alpha_6$ as indicated in the following Dynkin diagram:
\[
\xymatrix{
 \alpha_1 \ar@{-}[r] 
 & \alpha_3 \ar@{-}[r]
 & \alpha_4 \ar@{-}[r]
 & \alpha_5 \ar@{-}[r]
 & \alpha_6 \\
 && \alpha_2 \ar@{-}[u] &&
}
\]

\begin{ex}
As mentioned above, the simply connected group $E_6(\bbC)$ has two irreducible representations of dimension $27$, these are dual to each other and in the Bourbaki notation their highest weights are~$\varpi_1$, $\varpi_6$. Next there is a unique irreducible representation of dimension $78$, which is the adjoint representation with highest weight $\varpi_2$. Possibly after replacing the embedding $S\subset A$ by its negative we have
\[
 \varepsilon_{\varpi_d} \;\simeq\; 
 \begin{cases}
 \;\delta_S & \textnormal{if $d=1$}, \\
 \;\delta_{-S} & \textnormal{if $d=6$}, \\
 \;\delta_\Theta & \textnormal{if $d=2$},
 \end{cases}
\]
see~\cite{KraemerThreefolds}. We do not know a geometric description in the cases $d=3,4,5$.
%
\end{ex}

Although here we know much less than for Jacobians of curves, this still suffices for our purpose:

\begin{lem} \label{lem:muthreefold}
Any $\lambda\in \sfX^+$ satisfies 
$\mu \preceq \lambda$ for some $\mu \in \{ \varpi_1, \varpi_2, \varpi_6\}$.
\end{lem}

{\em Proof.} Any dominant weight can be written uniquely as a nonnegative linear combination 
\[
 \lambda \;=\; a_1 \varpi_1 + \cdots + a_6 \varpi_6
 \quad \textnormal{with} \quad a_1, \dots, a_6 \;\in\; \bbN_0.
\] 
Writing the fundamental weights as rational linear combinations of $\alpha_1, \dots, \alpha_6$ as in~\cite{BourbakiLie}, one has
$
\varpi_3 \succeq \varpi_6$,
$\varpi_4 \succeq \varpi_2$ and
$\varpi_5 \succeq \varpi_1$.
We successively apply these relations until we arrive at 
\[
 \nu \;=\; b_1 \varpi_1 + b_2\varpi_2 + b_6\varpi_6 \;\preceq\; \lambda
 \quad 
 \textnormal{with} \quad 
 \begin{cases}
 b_1 \;=\; a_1 + a_5, \\
 b_2 \;=\; a_2 + a_4, \\
 b_6 \;=\; a_3 + a_6.
 \end{cases}
\]
Using that $\varpi_2\succeq 0$ and $\varpi_1 + \varpi_6\succeq 0$, we can then further reduce to the dominant weight
\[ 
 \mu \;=\; c_1 \varpi_1 + c_2 \varpi_2 + c_6 \varpi_6 
 \;\preceq\; \nu 
 \quad \textnormal{where} \quad 
 (c_1, c_2, c_3) \;=\;
 \begin{cases}
 (0,1,0)\\
 \textnormal{or} \; (k,0,0)\\
 \textnormal{or} \; (0,0,k)
 \end{cases}
\] 
for some $k\in \bbN$. In the last two situations we can further reduce to the case $k=1$ by the dominance relations $2\varpi_{1}\succeq \varpi_6$, $3\varpi_1\succeq \varpi_2$, $2\varpi_6\succeq \varpi_1$ and $3\varpi_6 \succeq \varpi_2$.
\qed

\subsection{Proof of the main theorem for intermediate Jacobians}

Let $\Theta = X+Y$ be a decomposition as a sum of two geometrically nondegenerate subvarieties of positive dimension. Theorem~\ref{thm:basicstep} says that there exist dominant weights $\lambda, \mu \in \sfX^+$ such that
\[
 \cc(\delta_{nX}) \;=\; [n]_* \cc(\varepsilon_\lambda)
 \quad \textnormal{and} \quad 
 \cc(\delta_{nY}) \;=\; [n]_* \cc(\varepsilon_\mu).
\]
Remark~\ref{rem:weights} and lemma~\ref{lem:muthreefold} then imply \smallskip 
\begin{itemize}
\item $\dim X \geq 2$, with equality only if $nX = \pm nS$, \smallskip
\item $\dim Y \geq 2$, with equality only if $nY = \pm nS$. \smallskip
\end{itemize}
In fact we must have equality in both cases because $\dim X + \dim Y = \dim \Theta = 4$, and taking preimages under the isogeny $[n]: A\rightarrow A$ we obtain by irreducibility that
$X = \pm S + x$ and $Y = \pm S + y$
for some torsion points $x,y\in A[n]$. Finally, a look at weights or a direct geometric argument shows that the two signs must be opposite to each other, so theorem~\ref{thm:JT} follows. \qed

\medskip 

{\em Acknowledgements.} This paper originated in a question of Stefan Schreieder.
I would like to thank him and Mihnea Popa for interesting comments, and the referee for a very careful proofreading of the manuscript. 

\medskip

\bibliographystyle{amsplain}
\bibliography{Bibliography}

\providecommand{\bysame}{\leavevmode\hbox to3em{\hrulefill}\thinspace}
\providecommand{\MR}{\relax\ifhmode\unskip\space\fi MR }
\providecommand{\MRhref}[2]{%
  \href{http://www.ams.org/mathscinet-getitem?mr=#1}{#2}
}
\providecommand{\href}[2]{#2}
\begin{thebibliography}{10}

\bibitem{BBD}
{Beilinson, A., Bernstein, J. and Deligne, P.}, \emph{{Faisceaux Pervers}},
  Ast{\'e}risque \textbf{100} (1982).

\bibitem{BourbakiLie}
{Bourbaki, N.}, \emph{{Groupes et alg{\`e}bres de Lie --- Chapitres 4 {\`a}
  6}}, {Springer Verlag}, 2009, Reprint of the first edition from 1968.

\bibitem{BtD}
{Br\"ocker, T. and tom Dieck, T.}, \emph{{Representations of Compact Lie
  Groups}}, {Graduate Texts in Math.}, vol.~{98}, {Springer}, {1985}.

\bibitem{CMPSGenericVanishing}
{Casalaina-Martin, S., Popa, M. and Schreieder, S.}, \emph{{Generic vanishing
  and minimal cohomology classes on abelian fivefolds}}, {J.~Alg.~Geom.}
  \textbf{{27}} ({2018}), {553--581}.

\bibitem{CGIntermediateJacobian}
{Clemens, C. H. and Griffiths, P. A.}, \emph{{The intermediate Jacobian of the
  cubic threefold}}, {Annals of Math.} \textbf{95} (1972), 281--356.

\bibitem{DCM}
{De Cataldo, M. and Migliorini, L.}, \emph{{The decomposition theorem, perverse
  sheaves and the topology of algebraic maps}}, {Bull. Amer. Math. Soc.}
  \textbf{{46}} ({2009}), {535--633}.

\bibitem{DebMinimal}
{Debarre, O.}, \emph{{Minimal cohomology classes and Jacobians}}, {J. Alg.
  Geom.} \textbf{{4}} ({1995}), {321--335}.

\bibitem{DebarreComplexTori}
\bysame, \emph{{Complex tori and abelian varieties}}, {SMF/AMS Texts and
  Monographs}, vol.~{11}, {2005}.

\bibitem{DM}
{Deligne, P. and Milne, J. S.}, \emph{{Tannakian categories}}, Hodge Cycles,
  Motives, and Shimura varieties, Lecture Notes in Math., vol. 900, Springer
  Verlag, 1982, pp.~\mbox{101--228}.

\bibitem{GaL}
{Gabber, O. and Loeser, F.}, \emph{{Faisceaux pervers $\ell$-adiques sur un
  tore}}, Duke Math. J. \textbf{83} (1996), 501--606.

\bibitem{GW}
{Goodman, R. and Wallach, N. R.}, \emph{{Symmetry, representations and
  invariants}}, {Springer Verlag}, 2009.

\bibitem{Heinloth}
{Heinloth, F.}, \emph{{A note on functional equations for zeta functions with
  values in Chow motives}}, {Ann. Inst. Fourier Grenoble} \textbf{{57}}
  ({2007}), {1927--45}.

\bibitem{KraemerMicrolocal}
{Kr\"amer, T.}, \emph{{Characteristic cycles and the microlocal geometry of the
  Gauss map I}}, {\url{arXiv:1604.02389}}.

\bibitem{KraemerMicrolocalII}
{Kr{\"a}mer, T.}, \emph{{Characteristic cycles and the microlocal geometry of
  the Gauss map II}}, {\url{arXiv:1807.01929}}.

\bibitem{KraemerSemiabelian}
\bysame, \emph{{Perverse sheaves on semiabelian varieties}}, {Rend. Semin. Mat.
  Univ. Padova} \textbf{{132}} ({2014}), {83--102}.

\bibitem{KraemerThreefolds}
{{Kr{\"a}mer, T.}}, \emph{{Cubic threefolds, Fano surfaces and the monodromy of
  the Gauss map}}, {Manuscripta Math.} \textbf{{149}} ({2016}), {303--314}.

\bibitem{KrWSmall}
{Kr{\"a}mer, T. and Weissauer, R.}, \emph{{Semisimple super Tannakian
  categories with a small tensor generator}}, {Pacific J. Math.} \textbf{{276}}
  ({2015}), {229--248}.

\bibitem{KrWVanishing}
\bysame, \emph{{Vanishing theorems for constructible sheaves on abelian
  varieties}}, {J. Alg. Geom.} \textbf{{24}} ({2015}), {531--568}.

\bibitem{MartensSpecial}
{Martens, H. H.}, \emph{{On the variety of special divisors on a curve}}, {J.
  Reine Angew.~Math.} \textbf{{227}} ({1967}), {111--120}.

\bibitem{RanSubvarieties}
{Ran, Z.}, \emph{{On subvarieties of abelian varieties}}, {Invent. Math.}
  \textbf{62}, {459--480}.

\bibitem{SchnellHolonomic}
{Schnell, C.}, \emph{{Holonomic $\mathscr{D}$-modules on abelian varieties}},
  {Publ. Math. Inst. Hautes {\'E}tudes Sci.} \textbf{{121}} ({2015}), {1--55}.

\bibitem{SchreiederCurveSummands}
{Schreieder, S.}, \emph{{Theta divisors with curve summands and the Schottky
  problem}}, {Math. Annalen} \textbf{{365}} ({2016}), {1017--1039}.

\bibitem{SchreiederDecomposable}
\bysame, \emph{{Decomposable theta divisors and generic vanishing}}, {Int.
  Math. Res. Notices} \textbf{{16}} ({2017}), {4984--5009}.

\bibitem{WeBN}
{Weissauer, R.}, \emph{{Brill-Noether sheaves}}, {\url{arXiv:math/0610923}}.

\bibitem{WeissauerVanishing}
\bysame, \emph{{Vanishing theorems for constructible sheaves on abelian
  varieties over finite fields}}, {Math. Annalen} \textbf{{365}} ({2016}),
  {559--578}.

\end{thebibliography}

\end{document}